\documentclass{lms}
\usepackage{epic}
\usepackage{amsmath}[1996/11/01]
\usepackage{amssymb,amsfonts,latexsym,epsfig,graphics}

\DeclareMathOperator{\Stab}{Stab}
\DeclareMathOperator{\Syl}{Syl}
\DeclareMathOperator{\ad}{ad}
\DeclareMathOperator{\sign}{sign}

\DeclareMathOperator{\Sym}{Sym}
\DeclareMathOperator{\Alt}{Alt}

\DeclareMathOperator{\Hom}{Hom}

\DeclareMathOperator{\End}{End}
\DeclareMathOperator{\GL}{GL}

\DeclareMathOperator{\wt}{wt}

\DeclareMathOperator{\Aut}{Aut}

\newcommand{\disj}{\stackrel{.}{\cup}}

\newcommand{\Z}{{\mathbb{Z}}}
\newcommand{\Q}{{\mathbb{Q}}}
\newcommand{\F}{{\mathbb{F}}}

\newtheorem{theorem}{Theorem}[section] 
\newtheorem{lemma}[theorem]{Lemma}     
\newtheorem{corollary}[theorem]{Corollary}
\newtheorem{proposition}[theorem]{Proposition}

\newnumbered{assertion}{Assertion}    
\newnumbered{conjecture}{Conjecture}  
\newnumbered{definition}{Definition}
\newnumbered{hypothesis}{Hypothesis}
\newnumbered{remark}{Remark}
\newnumbered{note}{Note}
\newnumbered{observation}{Observation}
\newnumbered{problem}{Problem}
\newnumbered{question}{Question}
\newnumbered{algorithm}{Algorithm}
\newnumbered{example}{Example}
\newunnumbered{notation}{Notation} 

\title[ Automorphisms of doubly-even self-dual binary codes.]
 { Automorphisms of doubly-even self-dual binary codes.} 

\author{Annika G\"unther and Gabriele Nebe}

\classno{94B05 (primary), 20G25, 11E95 (secondary)}

\begin{document}
\maketitle

\begin{abstract}
The automorphism group of a binary
doubly-even self-dual code is always contained in the alternating group.
On the other hand, given a permutation group $G$ of degree $n$ there exists a doubly-even self-dual $G$-invariant code
if and only if $n$ is a multiple of $8$, every simple self-dual $\F_2G$-module occurs with
even multiplicity in $\F_2^n$, and
$G$ is contained in the alternating group.
\end{abstract}

\section{Introduction.}

Self-dual binary codes have become of great interest, also because of 
Gleason's theorem \cite{Gleason} that establishes a connection between 
coding theory and invariant theory of finite groups. 
Optimal self-dual codes often have the additional property of being 
doubly-even, which means that the weight of every codeword is divisible by 4 
(see Definition \ref{codes}).
It follows from Gleason's theorem that the length $n$ of a doubly-even
self-dual code $C\leq \F_2^n$ 
is a multiple of 8, see \cite[Theorem 3c]{MacWilliams-Sloane}, for instance.

This note studies the automorphism group
$
  \Aut(C):=\{ \pi \in \Sym_n \;|\; C\pi =C\}
$
of such a code.

Theorem \ref{codesalt} shows that the automorphism group of any
 doubly-even self-dual code is always 
contained in the alternating group, a very basic result which astonishingly
does not seem to be known.
On the other hand Theorem \ref{main_theorem} characterizes the permutation groups $G \leq \Sym _n$
that fix a doubly-even self-dual binary code. 
This result generalizes results by Sloane and Thompson \cite{ThompsSloane} and
Mart\'inez-P\'erez and Willems \cite{Willems}.

The first section considers codes as modules for their automorphism group.
The main result is the characterization of permutation groups that 
act on a self-dual code in Theorem \ref{codesWitt}. 
Section \ref{2adic} treats permutation groups as subgroups of the
2-adic orthogonal groups. The most important observation is Lemma 
\ref{sign} that expresses the sign of a permutation as a certain
spinor norm.
Given a self-dual doubly-even binary code $C$, the automorphism 
group of the even unimodular $\Z _2$-lattice obtained from $C$
by construction A (see Section \ref{A}) is contained in 
the kernel of this spinor norm. This immediately yields 
Theorem \ref{codesalt}. Theorem \ref{main_theorem} follows from 
this result together with Theorem \ref{codesWitt}.

\section{Codes.}

\begin{definition}\label{codes}
A {\em binary code} $C $ of length $n$
 is a linear subspace of $\F_2^n$.
 Let $b: \F_2^n \times \F_2^n \rightarrow \F_2, b(x,y):=\sum_{i=1}^nx_iy_i$
 be the standard scalar product. The {\em dual code} is
      $$
        C^{\perp}:=\{ v \in \F_2^n \;|\; b(v,c) = 0 \text{ for all } c \in C\}.
      $$
The code $C$ is called {\em self-orthogonal} if $C \subseteq C^{\perp}$ and {\em self-dual} if $C=C^{\perp}$.
      The {\em weight} $\wt(c)$ of a codeword $c \in C$ is the number of its nonzero entries.
      The code $C$ is called {\em doubly-even}, or {\em Type II}, if the weight of every word in $C$ is a multiple of 4.
\end{definition}

This section investigates binary linear codes as modules for a 
subgroup $G$ of their automorphism groups. 
The main result is Theorem \ref{codesWitt} that characterizes the permutation groups 
acting on some self-dual code. 
To this aim we need the representation theoretic notion of self-dual
modules, cf. Definition \ref{self-dual_modules}.
 Note that this paper uses two different notions of duality.
The dual of an $\F_2 G$-module $S$ over the finite group $G$ is the 
$\F_2G$-module $S^{\ast}=\Hom_{\F_2}(S,\F_2)$,
whereas the dual of a code $C \le \F_2^n$ is as in Definition \ref{codes}. 
For $G \le \Aut(C)$ the code $C$ is also an $\F_2G$-module,
which is represented with respect to a distinguished basis.

\begin{definition}\label{self-dual_modules}
      Let $S$ be a right $G$-module. Then the {\em dual module} $S^{\ast} = \Hom_{\F_2}(S,\F_2)$
      is a right $G$-module via $fg(s):=f(s g^{-1})$, for $f \in S^{\ast},\;g \in G$ and $s \in S$.
      If $S \cong S^\ast$ then $S$ is called {\em self-dual}.
\end{definition}

 \begin{theorem}\label{codesWitt}
        Let $G \le \Sym_n$. Then there exists a self-dual code $C \le \F_2^n$
        with $G \le \Aut(C)$ if and only if every self-dual simple $\F_2 G$-module $S$
        occurs in the $\F_2 G$-module $\F_2^n$ with even multiplicity.
  \end{theorem}

The proof of this theorem is prepared in a few lemmas.

\begin{lemma}\label{einzigchar2}
      Let $S$ be a simple self-dual $\F_2 G$-module, and assume that $S$ carries a non-degenerate symmetric $G$-invariant bilinear form
      $\varphi: S \times S \rightarrow \F_2$.
      Then $\varphi$ is unique, up to isometry.
\end{lemma}

\begin{proof}
Since $\varphi$ is non-degenerate and $G$-invariant, 
it induces an $\F_2G$-isomorphism $\alpha_{\varphi} : S \rightarrow S^{\ast},\; s \mapsto (s' \mapsto \varphi(s,s'))$.
     Let $\psi: S \times S \rightarrow \F_2$ be another non-degenerate symmetric $G$-invariant bilinear form on $S$, then
     $\alpha_{\psi}=\alpha_{\varphi}\circ \vartheta$ for some $\vartheta$ in the finite field $\mathfrak{E}:=\End_G(S)$ of all $\F_2G$-endomorphisms of $S$,
 and hence
     $$
       \psi(s,s')=\alpha_{\psi}(s)(s')=\alpha_{\varphi}(\vartheta(s))(s') = \varphi(\vartheta(s),s')
     $$
     for all $s,s' \in S$.
     Consider the involution $^{\ad}$ on $\mathfrak{E}$ given by $\varphi(s,\alpha(s'))=\varphi(\alpha^{\ad}(s),s')$, for $s,s' \in S$.
     Since both $\varphi$ and $\psi$ are symmetric we have
     $$
       \varphi(\vartheta(s),s')=\psi(s,s')=\psi(s',s)=\varphi(\vartheta(s'),s)=
\varphi(s,\vartheta(s')) = \varphi(\vartheta^{\ad}(s),s')
     $$
     for all $s,s' \in S$ and hence $\vartheta \in \mathfrak{F}=\{ \alpha \in \mathfrak{E} \;|\; \alpha^{\ad}=\alpha\}$.
 The involution $^{\ad}$ is either the identity on $\mathfrak{E}$ or a
field automorphism of order 2.
In the first case 
$\mathfrak{F} = \mathfrak{E} =\{ \alpha \alpha ^{\ad } = \alpha ^2 \mid \alpha \in 
\mathfrak{E} \}$ since squaring is an automorphism of the 
finite field $\mathfrak{E}$.
In the second case the map
$\mathfrak{E} \rightarrow \mathfrak{F},\; \alpha \mapsto \alpha\;\alpha^{\ad}$
is the norm map onto the fixed field $\mathfrak{F}$.
Hence in either case there exists some
 $\gamma \in \mathfrak{E}$ with $\gamma \gamma^{\ad}=\vartheta$. 
Now $\gamma$ induces an isometry
 between the spaces $(S,\varphi)$ and $(S,\psi)$ since
     $
      \psi(s,s')=\varphi( \vartheta(s),s')=\varphi(\gamma^{\ad }( \gamma(s)),s')=\varphi(\gamma (s),\gamma (s') )
     $
     for all $s,s' \in S$.
\end{proof}

\begin{lemma}\label{duality_orthogonality}
      Let $G \le \Sym_n$ and let $N \le M \le \F_2^n$ be $G$-submodules (i.e. $G$-invariant codes).
      Then $(M/N)^{\ast} \cong N^{\perp}/M^{\perp}$.
\end{lemma}

\begin{proof}
      Let $M_N^{\ast}:=\{f \in \Hom_{\F_2}(M,\F_2) \;|\; f(n) =0 \text{ for all } n \in N\} \le M^{\ast}$.
      Then $M_N^{\ast}$ is canonically isomorphic to $(M/N)^{\ast}$. Let
      $\beta : N^{\perp} \rightarrow M_N^{\ast},\;\; n' \mapsto (m \mapsto b(m,n'))$. Then $\beta $
      is well-defined and surjective, since $\Upsilon: \F _2^n \rightarrow M^{\ast},\;\; v \mapsto  (m \mapsto b(m,v))$
      is surjective, and $\Upsilon(v) \in M_N^{\ast}$ if and only if $v \in N^{\perp}$.
      Clearly $\beta $ has kernel $M^{\perp}$ and hence
      $N^{\perp}/M^{\perp} \cong M_N^{\ast} \cong (M/N)^{\ast}$.
\end{proof}

\begin{corollary}\label{metabolic_implies_even}
      Let $G \le \Sym_n$.
      If there exists a self-dual code $C \le \F_2^n$ with $G \le \Aut(C)$ then every self-dual simple $G$-module occurs
  with even multiplicity in
      a composition series of the $\F_2 G$-module $\F_2^n$.
\end{corollary}

\begin{proof}
      Let
      $
        C=N_k \ge N_{k-1} \ge \ldots \ge N_1 \ge N_0 =\{0\}
      $
      be a composition series of the $\F_2G$-module $C$. Then
      $$
       C=C^{\perp} = N_k^{\perp} \le N_{k-1}^{\perp} \le \ldots \le N_1^{\perp} \le N_0^{\perp } =\F_2^n
      $$
      is a composition series of $\F_2^n/C^{\perp}$, since
   dualizing 
yields an antiautomorphism $W\mapsto W^{\perp} $ of the submodule lattice
      of $\F_2^n$. The composition factors satisfy
      $$N_{i-1}^{\perp}/N_i^{\perp} \cong (N_i/N_{i-1})^{\ast},$$
      cf. Lemma \ref{duality_orthogonality}. Hence the claim follows.
\end{proof}

\begin{lemma}\label{order_two}
      Let $S$ be a simple self-dual $\F_2 G$-module endowed with a non-degenerate $G$-invariant symmetric bilinear form $\varphi$.
      The module $(U,\psi):=\perp_{i=1}^k (S,\varphi)$ contains a submodule $X$ with
      $$
        X=X^{\perp,\psi}:=\{u \in U \;|\; \psi(u,x)=0 \text{ for all } x \in X\}
      $$
      if and only if $k$ is even.
\end{lemma}

\begin{proof}
If $U$ contains such a submodule $X=X^{\perp,\psi}$
 then $k$ is even according to Corollary \ref{metabolic_implies_even}. 
Conversely,
if $k$ is even then $X:=\{(s_1,s_1,s_2,s_2,\ldots,s_{k/2},s_{k/2})\} \le U$ 
 satisfies $ X=X^{\perp,\psi}$.
\end{proof}


\begin{proof} (of Theorem \ref{codesWitt})
      If $C \le \F_2^n=:V$ is a self-dual $G$-invariant code then every self-dual simple module occurs 
with even multiplicity in a composition series of $V$
(see Corollary \ref{metabolic_implies_even}). Conversely, assume that every self-dual composition factor
      occurs in $V$ with even multiplicity, and
      let $M \le M^{\perp} \le V$ be a maximally self-orthogonal
 $G$-invariant code, i.e. there is no self-orthogonal $G$-invariant code in $V$ which properly
      contains $M$.

      On the $G$-module $M^{\perp}/M$ there exists a $G$-invariant non-degenerate symmetric bilinear form
      $$
        \varphi: M^{\perp}/M \times M^{\perp}/M \rightarrow \F_2,\;\; (m'+M,m''+M) \mapsto (m',m'').
      $$
      Any proper $\F_2G$-submodule $X$ 
of $(M^{\perp}/M,\varphi)$ with $X \subseteq X^{\perp,\varphi}$ (cf. Lemma \ref{order_two})
      would lift to a self-orthogonal $G$-invariant code in $V$ properly containing $M$, which we excluded in our assumptions.
      This implies that every $\F_2G$-submodule $X\le M^{\perp}/M$ has a $G$-invariant
complement $X^{\perp,\varphi}$, i.e. $M^{\perp}/M$ is
isomorphic to a direct sum of simple self-dual modules (see for instance \cite[Proposition (3.12)]{CurtisReiner}),
$(M^{\perp}/M,\varphi) \cong \perp_{S \cong S^{\ast}} (S,\varphi_S)^{n_S}$,
where $\varphi_S$ is a non-degenerate $G$-invariant bilinear form on $S$,
which is unique up to isometry by Lemma \ref{einzigchar2}.

      According to our assumptions, every simple self-dual $G$-module occurs with even multiplicity in $M^{\perp}/M$,
      i.e. all the $n_S$ are even. But this means that the $n_S$ must all be zero, according to
      Lemma \ref{order_two}, that is, $M=M^\perp$ is a self-dual code in $V$.
\end{proof}

The criterion in Theorem \ref{codesWitt} is not so easily tested. 
The next result gives a group theoretic condition that is sufficient for the existence 
of a self-dual $G$-invariant code. 
To this aim let $G \leq \Sym_n $ be a permutation group and write
$$\{ 1,\ldots , n \} = B_1 \disj \ldots \disj B_s $$ 
as a disjoint union of $G$-orbits and let $H_i := \Stab _G (x_i) $ be the
stabilizer in $G$ of some element $x_i \in B_i$ ($i=1,\ldots , s$).
For $1\leq i\leq s$ let $$m_i:=|\{ j\in \{1, \ldots , s\} \mid 
H_i \mbox{ is conjugate to } H_j \} | \mbox{ and  }
n_i:= [N_G(H_i) : H_i ].$$ 

\begin{proposition}\label{normstab} 
Assume that the product $n_i m_i $ is
even for all $1\leq i \leq s$. 
Then there is a $G$-invariant self-dual binary code $C \leq \F_2^n$.
\end{proposition}

\begin{proof}
If $H_i$ and $H_j$ are conjugate for some $i\neq j$
then the permutation representations 
of $G$ on $B_i$ and $B_j$ are equivalent and by Theorem \ref{codesWitt}
there is a self-dual $G$-invariant code in the direct sum 
$\F_2 ^{B_i \cup B_j} \cong \F_2^{|B_i|} \perp \F_2^{|B_j|}$ of two isomorphic $\F_2G$-modules.
It is hence enough to show the proposition for a transitive permutation 
group  $G \leq \Sym_n$ with stabilizer $H := \Stab_G(1) $ for which 
$[N_G(H):H] \in 2\Z$.
Let $(f_1,\ldots , f_n)$ be the standard basis of $\F_2^n$ such that 
$\pi \in \Sym _n$ 
maps $f_j $ to $f_{j \pi}$ for all $j = 1,\ldots,n$ and choose
$\eta  \in N_G(H) - H$  such that $\eta ^2 \in H$.
Put $N:= \langle H , \eta \rangle $ and 
$$G =  \disj _{s \in S} N s = 
  \disj _{s \in S} (H s \disj H \eta  s ).$$
Define 
$C:= \langle  f_{1s} + f_{1 \eta s }  : s \in S \rangle _{\F_2}.$
Then $C$ is a $G$-invariant code in $\F_2^n$ and  $C=C^{\perp }$ since
the given  basis of $C$ consists of 
$|S| = n/2$ pairwise orthogonal vectors of weight 2. 
\end{proof}

\section{From codes to lattices.}\label{A}

There is a well-known construction, called construction A  (see \cite[Section (7.2)]{ConwaySloane})
that associates to a pair $(R,C)$ of a ring $R$ with prime ideal $\wp $ and residue field
$R/\wp \cong \F $ and a code $C\leq \F ^n$ an $n$-dimensional lattice over $R$. 
We will apply this construction for binary codes and two different base rings:
$R = \Z$ and $R =\Z _2$, the ring of $2$-adic integers, where the prime ideal $\wp = 2R$ 
in both cases. 
So let $R$ be one of these two rings and let $K$  denote the field of fractions of $R$ and let 
 $V:=\langle b_1,\ldots , b_n \rangle _{K}$ be a vector space over $K$ 
with bilinear form defined by 
$$(\phantom{x},\phantom{x}) : V\times V \to K, (b_i,b_j) := \frac{1}{2} \delta _{ij}
= \left\{ \begin{array}{ll} 1/2 & i = j \\ 0 & i \neq j \end{array} \right. $$
and associated quadratic form $q:V\to K ,q(v) := \frac{1}{2} (v,v)$.
The orthogonal group of $V$ is 
$$O(V):= \{ g \in \GL(V) \mid (vg,wg) = (v,w) \mbox{ for all } v,w \in V \}.$$

\begin{definition}
A {\em lattice} $L\leq V$ is the $R$-span of a basis  of $V$.
The {\em dual lattice}
$$L^{\#} := \{ v\in V \mid (v,\ell ) \in R \mbox{ for all } \ell  \in L \} $$
is again a lattice in $V$.
$L$ is called {\em integral} if $L\subseteq L^{\#} $
or equivalently $(\ell _1,\ell _2) \in R$ for all $\ell _1,\ell _2 \in L$.
 $L$ is called {\em even} if $q(\ell ) \in R$ for all $\ell \in L$ and
{\em odd} if $L$ is integral and there is some $\ell \in L$ with $q(\ell ) \not\in R$.
$L$ is called {\em unimodular} if $L=L^{\#} $.
The {\em orthogonal group} of $L$ is
$$O(L) := \{ g\in O(V) \mid Lg = L \} .$$
\end{definition}

The following remark lists elementary properties of the lattice obtained from 
a code by construction A which can be seen by straightforward calculations.

\begin{remark}\label{remA}
Let  $M= \langle b_1,\ldots , b_n \rangle _{R} $ be the lattice  generated by the basis 
above and let $C \leq \F_2^n $ be a binary code.
Then the $R$-lattice
$$L:=A(R,C):= \{ \sum _{i=1}^n a_i b_i \mid a_i \in R ,
(a_1 +2R,\ldots ,a_n +2R)  \in C \}$$ 
is called the {\em codelattice} of $C$. 
Note that $2M \subset L \subset M$ and $L$ 
is the full preimage of $C \cong L/2M$ under the
natural epimorphism $M \to (R/2R)^n  = \F_2^n$. 
The lattice $L$ is even if and only if the code $C$ is doubly-even.
The dual lattice is $A(R,C)^{\#}=A(R,C^{\perp })$ and hence 
$L$ is unimodular if and only if $C$ is self-dual, and 
$L$ is an even unimodular lattice if and only if $C$ is a 
doubly-even self-dual code. 

The symmetric group $\Sym _n$ acts as orthogonal transformations on $V$ by
permuting the basis vectors. 
This yields an injective homomorphism 
$$\iota : \Sym _n \to O(V) ,\;\; \iota (\pi ): b_i \mapsto b_{i \pi} .$$
If $G = \Aut(C) $ is the automorphism group of $C$ then 
$\iota (G) \leq O(A(R,C))$.
\end{remark}

\section{Permutations as elements of the orthogonal group.}\label{2adic}

Let $\Q_2$ denote the field of $2$-adic numbers, 
$v_2: \Q _2 \to \Z \cup \{ \infty \} $ its natural valuation and 
$\Z _2 := \{ x\in \Q_2 \mid v_2(x) \geq 0 \}$ the ring of 
2-adic integers with unit group 
$\Z_2^* := \{ x\in \Z_2 \mid v_2(x) = 0 \}$. 
Let $V:=\langle b_1,\ldots , b_n \rangle _{\Q _2}$ 
be a bilinear space over $\Q _2$ of dimension $n>1$ as in Section \ref{A},
in particular $(b_i,b_j) = \frac{1}{2} \delta _{ij} $.
The orthogonal group 
$O(V)$
is generated by all reflections 
$$\sigma _v : V\to V , x\mapsto x - \frac{(x,v)}{q(v)} v $$
along vectors $v\in V$ with $q(v) \neq 0$ 
(see \cite[Satz (3.5)]{Kneser}, \cite[Theorem 43:3]{OMeara}).
Then the spinor norm defines a group homomorphism $h:O(V) \to C_2$ 
as follows:

\begin{definition}\label{defh}
Let 
$h:O(V) \to  C_2 = \{ 1, -1 \}$ be defined by 
$h(\sigma _v) := (-1) ^{v_2(q(v))} $ for all reflections $\sigma _v\in O(V)$.
Let $O^h(V) :=\{ g\in O(V) \mid h(g) = 1 \}$ denote the 
kernel of this epimorphism. 
\end{definition}

Note that the definition of $h$ depends on the chosen scaling of 
the quadratic form. 
It follows from the definition of the spinor norm 
(see \cite[Section 55]{OMeara}) that

\begin{lemma}
The map $h$ is a well-defined group epimorphism.
\end{lemma}

The crucial observation that yields the connection to coding
theory in Section \ref{main} is the following easy lemma.

\begin{lemma}\label{sign}
Let $\iota : \Sym_n \to O(V)$ be the homomorphism from Remark \ref{remA}.
Then $h\circ \iota = \sign $.
\end{lemma}

\begin{proof}
The symmetric group $\Sym_n$ is generated by transpositions $\tau _{i,j} = (i,j)$
for $i\neq j$.
Such a transposition interchanges $b_i$ and $b_j$ and fixes all 
other basis vectors and hence 
$\iota (\tau _{i,j} ) = \sigma _{b_i-b_j}$.
Clearly
$$h(\sigma _{b_i-b_j}) = (-1)^{v_2(q(b_i)+q(b_j))} = (-1)^{-1}= -1 = \sign(\tau _{i,j}).$$
\end{proof}


\begin{lemma} \label{ger}
Let $L\leq V$ be an even unimodular lattice.
Then $O(L) \leq O^h(V)$.
\end{lemma}

\begin{proof}
By \cite[Satz 4.6]{Kneser} the orthogonal group $O(L)$ is generated 
by reflections 
$$ O(L) = \langle \sigma _\ell  \mid \ell \in L, v_2(q(\ell )) = 0 \rangle .$$
Since $h(\sigma _\ell ) = (-1) ^{v_2(q(\ell ))} = 1$ for those vectors 
$\ell $, the result follows.
\end{proof}

We now assume that $n$ is a multiple of  8 and choose 
an orthonormal basis $(e_1,\ldots , e_n)$  of $V$ 
(i.e. $(e_i,e_j) = \delta _{ij}$).
Let $L:= \langle e_1,\ldots ,e_n \rangle _{\Z _2}$ be the 
unimodular lattice generated by these vectors $e_i$ and let
$$L_0 := \{ \ell \in L \mid q(\ell ) \in \Z_2 \}
=
\langle e_1+e_2,\ldots , e_1+e_n, 2 e_1 \rangle 
$$
be its even sublattice. 
Then 
$L_0^{\#} = \langle e_1,\ldots , e_{n-1},v:=\frac{1}{2} \sum _{i=1}^n e_i \rangle $.
Since $n$ is a multiple of 8 the vector
$2v\in L_0$ and $(v,v) = \frac{n}{4} $ is even.
Hence 
$L_0^{\#} /L_0\cong \F_2^2$ and the three 
lattices $L_i$ with $L_0 < L_i < L_0^{\#} $ corresponding to
the three 1-dimensional subspaces of $L_0^{\#}/L_0$ are given by 
$$L_1 := \langle L_0,v \rangle, \ L_2 := \langle L_0 , v-e_1 \rangle, \ 
L_3 = L .$$
Note that $L_1$ and $L_2$ are even unimodular lattices, whereas 
$L_3$ is odd.
In particular $O(L) = O(L_0)$ acts as the 
subgroup $$\{ 1,-1\} = C_2 \cong \{ I_2, \left( \begin{array}{cc} 1 & 1 \\ 0 & 1 \end{array} \right) \} \leq \GL_2(\F_2 )$$ on $L_0^{\#} / L_0 $ 
(with respect to the basis $(v+L_0,e_1+L_0)$). 
Let $f: O(L) \to C_2 = \{ \pm 1\} $ denote the resulting epimorphism.
So the elements in the kernel of $f$ (which equals $O(L)\cap O^h(V)$ as shown in 
the next lemma) fix both lattices $L_1$ and $L_2$ and all 
other elements in $O(L)$ interchange $L_1$ and $L_2$.

\begin{lemma} \label{unger} 
$f = h _{|O(L)} $
\end{lemma}

\begin{proof}
Let 
$R(L_0):= \langle \sigma _{\ell } \mid \ell \in L_0 , q(\ell ) \in \Z_2^*  \rangle $
be the reflection subgroup of $O(L_0)$.
By \cite[Satz 6]{Kneserpaper} $R(L_0) $ is the kernel of $f$. 
Since $h(\sigma _{\ell }) = 1$ for all $\sigma_{\ell } \in R(L_0)$, 
the group $R(L_0) \subset O(L) \cap O^h(V) $ is also
contained in the kernel of $h$.
The reflection $\sigma _{e_1}$ along the vector $e_1 \in L$ 
is in the orthogonal group $O(L) = O(L_0)$, interchanges the two 
lattices $L_1$ and $L_2$, and satisfies 
$h(\sigma _{e_1}) = -1$. 
Since $R(L_0)$ is a normal subgroup of index at most 2 in $O(L)$,
we obtain $O(L) = \langle R(L_0),\sigma _{e_1} \rangle $ 
and the lemma follows.
\end{proof}

\section{The main results.}\label{main}

\begin{theorem}\label{codesalt}
Let $C=C^{\perp}\leq \F_2^n$ be a doubly-even self-dual code.
Then  the automorphism group of $C$ is contained in the alternating group.
\end{theorem}

\begin{proof}
We apply construction A from Section \ref{A} to the code $C$ to obtain 
the codelattice $L:=A(\Z_2,C)$. 
By Remark \ref{remA} the lattice $L$ 
is  an even unimodular lattice. Hence by
Lemma \ref{ger} its orthogonal group $O(L) \leq O^{h}(V) $ is in the kernel of 
the epimorphism $h$ from Definition \ref{defh}.
 The image of $\Aut(C)$ under the homomorphism $\iota $ from 
Remark \ref{remA} is contained in $O(L)$,  hence
$\iota (\Aut(C)) \leq O(L) \leq O^h(V)$.
Since $h\circ \iota = \sign $ by Lemma \ref{sign} 
we have $\sign (\Aut(C)) = \{ 1 \}$ and therefore 
$\Aut (C) \leq \Alt _n $.
\end{proof}

\begin{theorem}\label{main_theorem}
Let $G\leq \Sym_n$. 
Then there is a self-dual doubly-even code $C=C^{\perp} \leq \F_2^n$ 
with $G\leq \Aut (C)$ if and only if the following three conditions 
are fulfilled:
\begin{itemize}
\item[(a)] $8\mid n$.
\item[(b)] Every self-dual composition factor of the $\F_2G$-module
$\F_2^n$ occurs with even multiplicity.
\item[(c)] $G\leq \Alt _n$.
\end{itemize}
\end{theorem}

\begin{proof}
\underline{$\Rightarrow:$} 
(a) is clear since the length of any doubly-even self-dual code
is a multiple of 8.
(b) follows from Theorem \ref{codesWitt} and 
(c) is a consequence of Theorem \ref{codesalt}. 
\\
\underline{$\Leftarrow :$} 
By Theorem \ref{codesWitt} the condition (b) implies the existence of
a self-dual code $X=X^{\perp}$ with $G\leq \Aut(X) $.
If $X$ is doubly-even then we are done. 
So assume that $X$ is not doubly-even and consider the codelattices
$$L := A(\Z,X) \mbox{ and } L_X := A(\Z _2,X) = L \otimes \Z _2 .$$
Then $L$ is a positive definite odd unimodular $\Z $-lattice 
 and hence its 2-adic
completion $L\otimes \Z_2 = L_X $ is an odd unimodular $\Z_2$-lattice
 having an orthonormal basis (see for instance \cite[Satz (26.7)]{Kneser}).
Hence $L_X$ is isometric to the lattice $L$  constructed 
just before Lemma \ref{unger}.
Since $G\leq \Alt _n$, 
the group $\iota (G) \leq O(L_X)$ lies in the kernel of the homomorphism
$f$ from Lemma \ref{unger} and therefore fixes the two even unimodular
lattices $L_1$ and $L_2$ intersecting $L_X$ in its even sublattice.
Let  $M= \langle b_1,\ldots , b_n \rangle _{\Z_2} $ be the lattice from Remark \ref{remA}
such that $2M < L_X < M $ and 
identify $M/2M = \bigoplus _{i=1}^n \Z_2/2\Z_2 b_i =  \bigoplus _{i=1}^n \F_2 b_i$ with 
$\F_2^n$. Then the
code $C:=L_1 / 2M  \leq  \F _2^n$ (such that $L_1=A(\Z_2,C)$) is 
a self-dual doubly-even code with $G\leq \Aut(C) $.
\end{proof}

\section{An application to group ring codes.}

As an application of our main Theorem \ref{main_theorem} we obtain a result (Theorem \ref{mainWil}) on the
existence of self-dual doubly-even binary group codes, given in \cite{ThompsSloane} and also in \cite{Willems}.
Binary group codes are ideals of the group ring $\F_2 G$, where $G$ is a finite group,
i.e. these are exactly the codes in $\F_2^{|G|}$ with $\rho_G(G) \le \Aut(C)$, where
$\rho_G : G \rightarrow \Sym_{G},\; g \mapsto (h \mapsto hg)$ is the regular representation of $G$.
Clearly $\rho_G(G) \leq \Alt _G $ if and only if the image $\rho_G(S)$ of any Sylow 2-subgroup $S \in \Syl_2(G)$ 
is contained in the alternating group.
Let $k:=[G:S]$ be the index of $S$ in $G$. Then $k$ is odd and the restriction of 
$\rho_G$ to $S$ is $(\rho_G )_{|S} = k \rho _S $. 
Hence $\rho_G(S) \leq \Alt _G $ if and only if $\rho_S(S) \leq \Alt _S $.

\begin{lemma}\label{Sylalt}
Let $S \neq 1 $ be a 2-group. Then $\rho_S(S) \leq \Alt _S $ if and only if $S$ is not cyclic.
\end{lemma}

\begin{proof}
If $S = \langle s \rangle$ is cyclic, then $\rho_S(s)$ is a $|S|$-cycle in $\Sym _S$ and hence its 
sign is -1 (because $|S|$ is even). 
On the other hand assume that $S$ is not cyclic. 
Then $S$ has a normal subgroup $N$ such that  $S/N \cong C_2 \times C_2$ is generated by elements
$aN,bN \in S/N$ of order $2$, with $abN=baN$.
Let $A := \langle a,N \rangle $ and $B = \langle b,N \rangle $.
Then $$S = \langle A,B \rangle  = A \disj bA = B \disj aB $$ 
and $b$ induces an isomorphism between the regular $A$-module $ A$ and $bA$,
so $A$ is in the kernel of the sign homomorphism.
Similarly $a$ gives an isomorphism between the regular $B$-module $B$ and $aB$,
so also $B$ is in the kernel of the sign homomorphism.
\end{proof}

The following observation follows from Proposition \ref{normstab} and is proven in \cite[Theorem 1.1]{Willems2}.

\begin{theorem} \label{Grpeven}
There is a self-dual binary group code $C\leq \F_2 G$ if and only if 
the order of $G$ is even.
\end{theorem}

\begin{proof}
\underline{$\Rightarrow $}: Clear, since $\dim (C) = \frac{|G|}{2}$  for any $C=C^{\perp } \leq \F_2G$.
\\
\underline{$\Leftarrow $}: Follows from  Proposition \ref{normstab}, because $\rho _G$ is a transitive permutation
representation and the the full group $G$ is the normalizer of the stabilizer $H:= \Stab_G(1) = {1}$.
\end{proof}

\begin{theorem}(\textup{see \cite{ThompsSloane},\cite{Willems}.})\label{mainWil}
      Let $G$ be a finite group. Then $\F_2 G$ contains a doubly-even self-dual group code if and only if the order of $G$ is
      divisible by $8$ and the Sylow $2$-subgroups  of $G$ are not cyclic.
\end{theorem}

\begin{proof}
The condition that the group order be divisible by 8 is equivalent to condition (a) of 
Theorem  \ref{main_theorem} and also implies (with Theorem \ref{Grpeven}) that there is some self-dual $G$-invariant code in 
$\F _2G$, which is equivalent to condition (b) of Theorem  \ref{main_theorem} by Theorem \ref{codesWitt}.
The condition on the Sylow 2-subgroups of $G$ is equivalent to 
$\rho_G(G) \leq \Alt_G$ by Lemma \ref{Sylalt} and hence
to condition (c) of Theorem  \ref{main_theorem}.
\end{proof}

Our last application concerns the automorphism group $G=\Aut(C) $ of a putative extremal Type II 
code $C$ of length 72. 
The paper \cite{XXX}  shows that 
any automorphism of $C$ of order 2 acts fixed point freely, so
any Sylow 2-subgroup $S$ of $G$ acts as a multiple of 
the regular representation.
In particular $|S|$ divides 8.
Our results show that $S$ is not cyclic of order 8, which already
follows from \cite[Theorem 1]{ThompsSloane}.

\begin{corollary}
Let $C$ be a self-dual doubly-even binary code of length 72 with minimum distance 16.
Then $C$ does not have an automorphism of order 8.
\end{corollary}

\section{A characteristic 2 proof of Theorems \ref{codesalt} and \ref{main_theorem}}

As remarked by Robert Griess 
 one may prove Theorem \ref{codesalt} and \ref{main_theorem} without using 
characteristic 0 theory. 

Assume that $n$ is a multiple of 8 and
 let ${\bf 1} := (1,\ldots , 1) \in \F _2^n $ denote the 
all ones vector. Then
$$V= {\bf 1}^{\perp }/\langle {\bf 1} \rangle = \{ x\in  \F_2 ^n \mid
\wt (x) \mbox{ is even } \}/\langle {\bf 1} \rangle $$
becomes a quadratic  module of dimension $n-2$ over $\F_2 $ by putting
$$q:V \to \F_2, \overline{x}:=x+\langle {\bf 1} \rangle \mapsto \frac{1}{2} \wt (x) + 2\Z .$$
The associated bilinear form 
$b(\overline{x},\overline{y}) = q(\overline{x}+\overline{y})-q(\overline{x})-q(\overline{y}) = x\cdot y $ is inherited from the standard inner product 
and the maximal isotropic subspaces of $V$ are the images of the
doubly-even self-dual codes in $\F_2^n$.

The orthogonal group $O(V,q) \cong O_{n-2}^+(2)$ acts transitively on the 
set of maximal isotropic subspaces of $V$.
Fix one such subspace $U$. 
Then the {\em Dickson invariant} is 
$$D:O(V,q) \to \{ 1,-1\} ; D(g) := (-1) ^{\dim(U/U\cap Ug)} $$
a well-defined homomorphism that does not depend on the choice of $U$
(\cite[Theorem 11.61]{Taylor}).

The symmetric group $\Sym _n$ acts by coordinate permutations on 
$\F_2^n$. Since ${\bf 1} \pi = {\bf 1} $ for all $\pi \in \Sym _n$ 
and permutations preserve the weight this gives rise to an embedding
$\iota : \Sym _n \to O(V,q)$.
The following lemma also follows from the geometric characterization
of the Dickson invariant in \cite[p. 160]{Taylor}
(see also \cite{Dye} and \cite{Dieudonne}).

\begin{lemma}
$D\circ \iota = \sign $.
\end{lemma}

\begin{proof}
It is enough to find a transposition that is not in the kernel of 
the Dickson invariant. To this aim choose the Type II code $C$ with 
generator matrix 
$$ \left( \begin{array}{ccccccccccc}
1 & 1 & 1 & 1 & 0 & 0 & 0 & \ldots & 0 & 0 & 0 \\
1 & 1 & 0 & 0 & 1 & 1 & 0 & \ldots & 0 & 0 & 0 \\
\vdots & \vdots & \vdots & \vdots &
\ddots & \ddots & \ddots & \ddots & \ddots & \vdots & \vdots \\
1 & 1 & 0 & 0 & 0 & 0 & 0 & \ldots & 0 & 1 & 1 \\
1 & 0 & 1 & 0 & 1 & 0 & \ldots & \ldots & \ldots & 1 & 0 
\end{array} \right) $$
and let $U:= C/\langle \bf 1 \rangle $.
Then $U \iota (\tau _{1,2}) \cap U$ has co-dimension 1 in U.
\end{proof}

Now we can use the Dickson invariant $D$ to replace the spinor norm 
$h$ to obtain the main results.
It is immediate that 
$\Stab _{O(V,q)} (U) \subset \ker (D) $ 
(see also \cite[Exercise 11.19]{Taylor}) from which one obtains 
Theorem \ref{codesalt}.

The proof of Theorem \ref{main_theorem} can also be modified.
Condition (b) implies the 
existence of a self-dual $G$-invariant 
code $X$. If $X$ is doubly-even, then we are done; if not, then 
let $X_0 := \{ x\in X \mid \wt (x) \in 4\Z \}$ denote the doubly-even subcode
of $X$. This is a subcode of codimension 1 in $X$ and 
$X_0^{\perp }/X_0 \cong \F_2\oplus \F_2 $ is of dimension 2. 
Since the length of $X$ is divisible by 8, the full preimages $C_1$ and $C_2$ of
the  other  two non-trivial
subspaces of $X_0^{\perp }/X_0$ both are self-dual doubly-even codes.
Since the co-dimension of the intersection $\dim (C_i/(C_1\cap C_2)) = 1$ 
is odd, any permutation $\pi $ with $C_1 \pi  = C_2$ has to have 
$\sign (\pi ) = D(\iota (\pi )) = -1 $. 
Since $G \leq \Alt _n$, all elements of $G$ have to fix both 
codes $C_1$ and $C_2$ and hence these yield $G$-invariant 
doubly-even self-dual codes.

The proof of Theorem \ref{codesalt} given here directly generalizes 
to generalized doubly-even codes as well as to odd characteristic.
Note that in odd characteristic the Dickson invariant is the same as
the determinant of an orthogonal mapping.
For further details we refer to the first author's thesis.

\begin{theorem}
(a) Let $C=C^{\perp } \leq \F_{2^d}^n$ be a generalized doubly even 
code as defined in \cite{Quebbemann}.
Then $P(C) \leq \Alt _n$.
\\
(b) Let $q$ be an odd prime power and $C=C^{\perp } = 
\{ x\in \F_q^n \mid \sum _{i=1}^n x_i c_i = 0 \mbox{ for all } 
c\in C \} $. 
Then any monomial automorphism $g\in \Stab _{C_2\wr S_n} (C) $
has determinant 1.
\end{theorem}

\affiliationone{
A. G\"unther and G. Nebe,  \\
Lehrstuhl D f\"ur Mathematik, \\ RWTH Aachen University
52056 Aachen, Germany
   \email{
annika.guenther@math.rwth-aachen.de, gabriele.nebe@math.rwth-aachen.de 
}}
\end{document}